 \documentclass[final,1p,times]{elsarticle}
\usepackage{amsmath,amssymb,url}
\usepackage{enumitem} 
\numberwithin{equation}{section}

\usepackage{amssymb}
\usepackage{amsthm}
\usepackage{hyperref}
\usepackage{amsmath,amssymb,amsopn,amsfonts,mathrsfs,amsbsy,amscd}
\usepackage{longtable}
\journal{Proceedings of the Edinburgh Mathematical Society}

\newtheorem{definition}{Definition}[section]
\newtheorem{theorem}{Theorem}[section]
\newtheorem{proposition}{Proposition}[section]

\numberwithin{equation}{section}

\begin{document}
\bibliographystyle{amsplain}
\begin{frontmatter}
\title{Local superderivation and super-biderivation on generalized quaternion algebra}

\author[label2]{Hassan Oubba}
\address[label2]{Moulay Ismail University of Meknès\\ Faculty of Sciences Meknès
B.P. 11201 Zitoune Meknes, Morocco\\
{hassan.oubba@edu.umi.ac.ma}}

\begin{abstract}
Let $\mathcal{H}^{a,b}$ be the generalized quaternion algebra over a unitary commutative ring. This paper aims  to investigate super-biderivations and local superderivations on the generalized quaternion algebra, which is viewed as a class of Lie superalgebra. It turns out that on generalized quaternion algebras, any local superderivation is a superderivation. 
\end{abstract}
\end{frontmatter}

\section{Introduction}
Throughout the paper, $\mathcal{R}$ will denote a commutative ring with unity and all algebras and modules will be unital over $\mathcal{R}$. Let $A$ be an algebra (not necessarily associative) and $Z(\mathcal{A})$ denote the center of $A$. As usual, the Lie product is denoted by $[x,y]=xy-yx$,  for all $x,y\in\mathcal{A}$. A linear map $d:\mathcal{A}\to\mathcal{A}$ is a derivation  if $d(xy)=d(x)y+xd(y)$  for all $x,y\in\mathcal{A}$. As generalizations of derivations, biderivations and local derivations occupy an essential position in
 the research of the structure theory of algebras, which have aroused a great many authors interests. The exploration of local and 2-local derivations of non-associative algebras was pioneered in works by Ayupov and Kudaybergenov, with specific attention to Lie algebras in \cite{Ayo1, Ayo2}. They demonstrated that there are no non-trivial local and 2-local derivations on semisimple finite-dimensional Lie algebras. Over recent years there has been considerable interest in the study of local derivation on algebra,  for various algebras, there
 have been a lot of investigations done to find conditions implying that a local derivation  is actually a
 derivation, or construct examples of local derivations which are not derivations (see \cite{Ayo1, Ayo2, brl, oubw, oubq, oubd, oubm}. The same questions were also considered in the context of the super case \cite{fan, xia}.\\
 The
 notion of biderivations was introduced to Lie algebra in \cite{wang}. For the last few years, there are a lot of
 results on biderivations of various algebras (cf. \cite{wang, brb, Tang, bn, bbo, oubm, dil, xia} and references therein). Furthermore,
 the concept of biderivations on Lie algebras is generalized to the super case in \cite{fan, xia} independently.
 Up to now, there are several efforts on determining the super-biderivations of certain Lie superalgebras
 (cf. \cite{fan, xia}).

 In \cite{chen} the authors have shown that  every local superderivation on basic classical Lie superalgebras except for $A(1,1)$ over the complex number field $\mathbb{C}$ is a superderivation. Also, in \cite{fan} the authors have prove that all super-biderivations on the centerless super-Virasoro algebras are inner super-biderivations. This paper considers the analogous problem of innerness of local superderivations and super-biderivations of generalized quaternion algebras as a class of Lie superalgebra. We first show that every local superderivation of $\mathcal{H}^{a,b}$ is  superderivation. Moreover, we give a characterization of super-biderivations of generalized quaternion algebra.\\

 This paper is organized as follows. In Section 2, we review the basic notations which are used
 in this paper. In section 3, we prove that any local superderivation on the generalized quaternion algebra is a derivation. In Section 4, we compute the super-biderivations on generalized quaternion algebra is a derivation.

\section{Preliminaries}
Throughout the paper, by an algebra, we shall mean an algebra over a fixed unital commutative ring $\mathcal{R}$. We assume without further mention $\frac{1}{2}\in \mathcal{R}$. A superalgebra is a $\mathbb{Z}_{2}$-graded algebra. This means that there exist $\mathcal{R}$-submodules $\mathcal{A}_{0}$ and $\mathcal{A}_{1}$ of $\mathcal{A}$ such that $\mathcal{A}=\mathcal{A}_{0}\oplus\mathcal{A}_{1}$ and $\mathcal{A}_{\alpha}\mathcal{A}_{\beta}\subseteq\mathcal{A}_{\alpha+\beta},$ where indices are computed modulo $2$. A superalgebra is called trivial if $\mathcal{A}_{1}=0$. Elements in $\mathcal{A}_{0}\cup\mathcal{A}_{1}$ are said to be homogeneous of degree and we write $|x|=\alpha$ to mean $x\in\mathcal{A}_{\alpha}$. We say that $\mathcal{A}_{0}$ is the even part of $\mathcal{A}$ and $\mathcal{A}_{1}$ is the odd part of $\mathcal{A}$. Define a product in $\mathcal{A}_{0}\cup\mathcal{A}_{1}$ the Lie superproduct, by
$$[x,y]_{s}=xy-(-1)^{|x||y|}yx$$
for $x,y\in\mathcal{A}_{0}\cup\mathcal{A}_{1}$. The Jordan superproduct is defined by
$$x\circ_{s}y=\frac{1}{2}(xy+(-1)^{|x||y|}yx)$$
Extend $[x,y]_{s}$ and $x\circ_{s}y$ by bilinearity to $\mathcal{A}\times\mathcal{A}$ accordingly, and
\begin{align*}
[x,y]_{s} &= [x_{0},y_{0}]_{s}+[x_{1},y_{0}]_{s}+[x_{0},y_{1}]_{s}+[x_{1},y_{1}]_{s} \\
x\circ_{s}y &= x_{0}\circ_{s}y_{0}+x_{1}\circ_{s}y_{0}+x_{0}\circ_{s}y_{1}+x_{1}\circ_{s}y_{1},
\end{align*}
where $x=x_{0}+x_{1},$ $y=y_{0}+y_{1}$. Note that in trivial superalgebras the Lie superproduct coincides with the Lie product and the Jordan superproduct coincides with the Jordan product.

Generalized quaternions and quaternion algebras have been introduced in the last decades as tools for studying quadratic form theory. This construction is essentially a natural generalization of $\mathcal{H}=\mathcal{H}(\mathbb{R})=\mathbb{R}\oplus\mathbb{R}i\oplus\mathbb{R}j\oplus\mathbb{R}k$, the well-known quaternion algebra over the real numbers introduced by sir Hamilton. In [\cite{lam}, Chapter III], Lam works with quaternion algebras with coefficients over an arbitrary field of characteristic distinct from 2. More generally, preserving Lam's notations, we will consider here quaternion algebras over a commutative ring with $\frac{1}{2}\in \mathcal{R}$.

Let $\mathbf{R}$ be a ring with unity and $a, b$ are units. The generalized quaternion algebra denoted by $\mathcal{H}^{a,b}$ consists of elements of the form $x+yi+zj+wk$ where $x, y, z, w\in \mathcal{R}$, with componentwise addition and multiplication obeying the relations
\begin{align*}
i^{2}&=-a, \quad j^{2}=-b, \quad k^{2}=-ab \\
ij&=-ji=k \\
jk&=-kj=b i \\
ki&=-ik=a j
\end{align*}
Clearly, $\mathcal{H}^{a,b}$ is a noncommutative associative algebra with unity and a two-sided module over $\mathcal{R}$ with basis $\{1,i,j,k\}$. Specifically, when $a=b=1,$ we denote the algebra by $\mathcal{H}(R)$ in reference to Hamilton's quaternions. By letting $\mathcal{A}_{0}=\mathcal{R}\oplus \mathcal{R}i$ and $\mathcal{A}_{1}=\mathcal{R}j\oplus \mathcal{R}k$, it is easily verified that $\mathcal{H}^{a,b}=\mathcal{A}_{0}\oplus\mathcal{A}_{1}$ becomes a superalgebra. In the following, we always consider the above-mentioned superalgebraic structure on a quaternion algebra.

The following definition will be needed throughout the paper. For $\alpha=0,1,$ recall that a superderivation of degree $\alpha$ is a linear map $D_{\alpha}: \mathcal{A}\to \mathcal{A}$ such that $D_{\alpha}(\mathcal{A}_{\beta})\subseteq\mathcal{A}_{\alpha+\beta}$ (index modulo 2) and
$$D_{\alpha}(xy)=D_{\iota}(x)y+(-1)^{\alpha|x|}xD_{\alpha}(y)$$
for all $x,y\in\mathcal{A}_{0}\cup\mathcal{A}_{1}$. A superderivation of $\mathcal{A}$ is the sum of a superderivation of degree $0$ and a superderivation of degree $1$. Note that every superderivation of degree $0$ is a derivation from $\mathcal{A}$ to $\mathcal{A}$. For example, for any $x_{0}\in\mathcal{A}_{0}$ and $x_{1}\in\mathcal{A}_{1}$, the maps $I_{x_{0}}(y)=[x_{0},y]_{s}$ and $I_{x_{1}}(y)=[x_{1},y]_{s}$ are superderivations of degree $0$ and $1$, respectively. Hence for any $x=x_{0}+x_{1}\in\mathcal{A}$ the linear mapping $I_{x}(y)=[x,y]_{s}$ is a superderivation, since $[x,y]_{s}=[x_{0},y]_{s}+[x_{1},y]_{s}$. The superderivation $I_{x}(y)=[x,y]_{s}$ will be called the inner superderivation. The $\mathcal{R}$-module of even (odd) superderivations on $\mathcal{A}$ is a $\mathcal{R}$-submodule of $\mathcal{L}(\mathcal{A})$ (the associative algebra of linear maps from $\mathcal{A}$ to $\mathcal{A}$). We denote the direct sum of these two $\mathcal{A}$-submodules by $\mathrm{Der}_{s}(\mathcal{A})$ and call it the $\mathcal{R}$-module of superderivations. The $\mathcal{R}$-module $\mathrm{Der}_{s}(\mathcal{A})$ with the following bracket is a Lie superalgebra:
$$[D,D^{\prime}]_{s}=DD^{\prime}-(-1)^{|D||D^{\prime}|}D^{\prime}D,$$
where $D, D^{\prime}$ are homogeneous superderivaions. On the other hand, the inner superderivations $\mathrm{Inn}_{s}(\mathcal{A})$ form a Lie ideal. The quotient superalgebra $\mathrm{Der}_{s}(\mathcal{A})/\mathrm{Inn}_{s}(\mathcal{A})$ is called the outer superderivation algebra of $\mathcal{A}$ and is denoted by $\mathrm{Out}_{s}(\mathcal{A})$.

 Let $x$ be an
element in $\mathcal{H}^{a,b}$. Then we can write
$$x=x_{1}1+x_{2}i+x_{3}j+x_{4}k$$
for some elements $x_{1}$, $x_{2}$, $x_{3}$, and $x_{4}$ in $\mathcal{R}$. Throughout the paper, let
$$\overline{x}=(x_{1},x_{2},x_{3},x_{4})^{T}.$$
Conversely, if $v=(v_{1},v_{2},v_{3},v_{4})^{T}$ is a column vector with $v_{1}$, $v_{2}$, $v_{3}$ and $v_{4}$ in $\mathcal{R}$
then, throughout the paper, we will denote by the element
$$v_{1}1+v_{2}i+v_{3}j+v_{4}k;$$
i.e.,
$$\widehat{v}=v_{1}1+v_{2}i+v_{3}j+v_{4}k.$$

\section{Local and 2-local superderivations on $\mathcal{H}^{a,b}$}
This section is devoted to the first main theorem of this paper. We prove that any local superderivation on the generalized quaternion algebra is a superderivation.\\

Our principal tool for the description of local superderivations and super-biderivations of $\mathcal{H}^{a,b}$ is the following two
propositions.
\begin{proposition}\cite{he}\label{pr1}
    Let $d$ be a superderivation of degree $0$ on the generalized quaternion algebra $\mathcal{H}^{a,b}$. Then the matrix representation $M$ of $d$ is as follows:
$$M=\begin{pmatrix}0&0&0&0\\ 0&0&0&0\\ 0&0&0&-a \lambda\\ 0&0&\lambda&0\end{pmatrix}$$
where $\lambda,\mu,\nu\in \mathcal{R}$. Here the action of $M$ corresponds to multiplying the matrix by a column on the right.
\end{proposition}
\begin{proposition}\cite{he}\label{pr2}
    Let $d$ be a superderivation of degree $1$ on the generalized quaternion algebra $\mathcal{H}^{a,b}$. Then the matrix representation $M$ of $D$ is as follows:
$$M=\begin{pmatrix}0&0&\mu&\nu\\ 0&0&0&0\\ 0&-b^{-1}\nu&0&0\\ 0&b^{-1}\mu&0&0\end{pmatrix}$$
where $\mu,\nu\in \mathcal{R}$. Here the action of $M$ corresponds to multiplying the matrix by a column on the right.
\end{proposition}
\begin{definition}
    A linear map $\Delta: \mathcal{H}^{a,b} \times \mathcal{H}^{a,b} \to \mathcal{H}^{a,b} $ is called local sperderivation of degree $\gamma \in \mathbb{Z}_2$ if for any $x\in \mathcal{H}^{a,b}$ there exists a derivation $d_x$ of $\mathcal{H}^{a,b}$ of degree $\gamma$ such that 
    $$\Delta(x)=d_u(x).$$
\end{definition}
\begin{theorem}\label{th1}
   Each local superderivation of degree $0$ on the  algebra $\mathcal{H}_{\alpha, \beta}$ is a superderivation of degree $0$. 
\end{theorem}
\begin{proof}
 Consider an arbitrary local superderivation $\Delta$ of $\mathcal{H}^{a,b}$ and express it as
$$\Delta(x)=M\overline{x}, \hspace{0.3 cm} x \in \mathcal{H}^{a,b},$$
where $M=(m_{ij})_{1 \leq i,j \leq 4}$ represents the matrix of $\Delta$ in the basis $\lbrace 1,i,j,k \rbrace$, and $\overline{x}=(x_1,x_2,x_3,x_4)$ is the vector corresponding to $x$. For any $x \in \mathcal{H}^{a,b}$, there exists an element $\lambda^x$ (depending on $x$) in $\mathcal{R}$ such that

$$M\overline{x}=\begin{pmatrix}
0&0&0&0\\
0&0 &0&0\\
0& 0& 0& -a \lambda^x \\
0& 0 & \lambda^x & 0
\end{pmatrix} 
\begin{pmatrix}
x_1\\
x_2\\
x_3\\
x_4
\end{pmatrix}$$
\begin{gather*}
		  	 \begin{cases}
		  	m_{11}x_1+m_{12}x_2+m_{13}x_3+m_{14}x_4=0 ;&\\
		  	\,m_{21}x_1+m_{22}x_2+m_{23}x_3+m_{24}x_4=0; \\
		  	\,m_{31}x_1+m_{32}x_2+m_{33}x_3+m_{34}x_4=-a\lambda^xx_4; &\\	
		  	\, m_{41}x_1+m_{42}x_2+m_{43}x_3+m_{44}x_4=\lambda^x x_3.&\\ 	  	 
		  	 \end{cases}
\end{gather*}
Taking $x=1=\widehat{(1,0,0,0)}$, we have
\begin{gather*}
		  	 \begin{cases}
		  	m_{11}=0 ;&\\
		  	\,m_{21}=0; \\
		  	\,m_{31}=0;&\\	
		  	\, m_{41}=0.&\\ 	  	 
		  	 \end{cases}
\end{gather*}
For $x=i=\widehat{(0,1,0,0)}$, we get
\begin{gather*}
		  	 \begin{cases}
		  	m_{12}=0 ;&\\
		  	\,m_{22}=0;&  \\
		  	\,m_{32}=0  ;&\\	
		  	\, m_{42}=0 .&\\ 	  	 
		  	 \end{cases}
\end{gather*}
Taking also $x=j=\widehat{(0,0,1,0)}$, we deduce
\begin{gather*}
		  	 \begin{cases}
		  	m_{13}=0 ;&\\
		  	\,m_{23}=0 ; \\
		  	\,m_{33}=0  ;&\\	
		  	\,m_{43}=\lambda^j .&\\ 	  	 
		  	 \end{cases}
\end{gather*}
Finally taking $x=k=\widehat{(0,0,0,1)}$, we get
\begin{gather*}
		  	 \begin{cases}
		  	m_{14}=0 ;&\\
		  	\,m_{24}=0 ; \\
		  	\,m_{34}=-a\lambda^k ;&\\	
		  	\, m_{44}=0 .&\\ 	  	 
		  	 \end{cases}
\end{gather*}
Therefore,
$$M=\begin{pmatrix}
    0&0&0&0\\
    0&0&0&0\\
    0&0&0&-a\lambda^j\\
    0&0&\lambda^k&0
\end{pmatrix}$$
Since, $\Delta$ is a linear map, then
\begin{equation}\label{eq1}
    \Delta(j+k)=\Delta(j)+\Delta(K),
\end{equation}
and 
$$\Delta(j+k)=\lambda^{j+k}k-a \lambda^{j+k}j,$$
By (\ref{eq1}), we have
$$\lambda^{j+k}=\lambda^{j}, \quad-a\lambda^{j+k}=-a\lambda^{k}.$$
This implies that $\lambda^j=\lambda^k$ and 
$$M=\begin{pmatrix}
    0&0&0&0\\
    0&0&0&0\\
    0&0&0&-a\lambda^j\\
    0&0&\lambda^j&0
\end{pmatrix}$$
Hence, by Proposition \ref{pr1}, $\delta$ is a superderivation. This completes the proof.
\begin{theorem}\label{th2}
   Each local superderivation of degree $1$ on the  algebra $\mathcal{H}_{\alpha, \beta}$ is a superderivation of degree $1$. 
\end{theorem}
\begin{proof}
    The proof of this theorem is similar to the proof of Theorem \ref{th1}.
\end{proof}
\section{super-biderivation on $\mathcal{H}^{a,b}$}
This section, is devoted to caracterezation of super-biderivation on generalized quaternion algebra.
\begin{definition}
 A bilinear map $\delta : \mathcal{H}^{a,b} \times \mathcal{H}^{a,b} \rightarrow \mathcal{H}^{a,b}$ is called a \emph{super-biderivation} of $\mathcal{H}^{a,b}$ if
\begin{align*}
&\delta(xy, z) = (-1)^{|\delta||x|}x \delta(y, z) + (-1)^{|y||z|}\delta(x, z) y; \\
&\delta(x, yz) = \delta(x, y), z + (-1)^{(|\delta|+|x|)|y|}y\delta(x, z);
\end{align*}
for all homogeneous $x, y, z \in \mathcal{H}^{a,b}$.
\end{definition}

A super-biderivation $\delta$ of $\mathcal{H}^{a,b}$ is of degree $\gamma \in \mathbb{Z}_2$ if $\delta$ is a super-biderivation such that $\delta({\mathcal{H}^{a,b}}_\alpha, {\mathcal{H}^{a,b}}_\beta) \subseteq {\mathcal{H}^{a,b}}_{\alpha+\beta+\gamma}$ for any $\alpha,\beta \in \mathbb{Z}_2$. Denote by $\mathrm{BDer}_\gamma(\mathcal{H}^{a,b})$ the set of all super-biderivations of degree $\gamma$ of $\mathcal{H}^{a,b}$. We have,
\[
\mathrm{BDer}(\mathcal{H}^{a,b}) = \mathrm{BDer}_0(\mathcal{H}^{a,b}) \oplus \mathrm{BDer}_1(\mathcal{H}^{a,b}),
\]
where $\mathrm{BDer}(\mathcal{H}^{a,b})$ represent the collection of all super-biderivations defined on $\mathcal{H}^{a,b}$.

Now, let $\delta$ be a super-biderivation of  $\mathcal{H}^{a,b}$ and $x, y \in \mathcal{H}^{a,b}$, such that $x =x_1+x_2i+x_3j+x_4k$ and $y = y_1+y_2i+y_3j+y_4k$. Then, by the bilinearity of $\delta$, we obtain,
\[
\delta(x, y) = x_1\delta(1,y)+x_2\delta(i,y)+x_3\delta(j,y)+x_4\delta(k,y).
\]
It is easy to see that for any $x\in \mathcal{H}^{a,b}$ the linear maps $\delta(x,.):\mathcal{H}^{a,b}\to \mathcal{H}^{a,b}, \, y\mapsto \delta(x,y)$ and $\delta(.,x):\mathcal{H}^{a,b}\to \mathcal{H}^{a,b}, \, y\mapsto \delta(y,x)$ are superderivations of $\mathcal{H}^{a,b}$.
\subsection{Characterization of $BDer_0(\mathcal{H}^{a,b})$}
Now, we state and prove our main theorem of  superderivations of degree $0$ on
generalized quaternion algebra.
\begin{theorem}\label{thb1}
    Let $\delta$ be a skew-super-biderivation of degree $0$ on $\mathcal{H}^{a,b}$. Then, there exists $\lambda \in \mathcal{R}$ such that 
    $$\delta(x,y)=-ab\lambda(x_3y_3+x_4y_4)+a\lambda(x_4y_2-x_2y_4)j+\lambda(x_2y_3-x_3y_2)k, \quad \forall x,y \in \mathcal{H}^{a,b}.$$
\end{theorem}
\begin{proof}
Let $\delta$ be a skew-super-biderivation of degree $0$ on  $(\mathcal{H}^{a,b})$, then for any $x,y \in \mathcal{H}^{a,b}$ we have 
\begin{eqnarray*}
    &&\delta(1,xy)=\delta(1,x)y+x\delta(1,y),\\
    &&\delta(i,xy)=\delta(i,x)y+x\delta(i,y),\\
    &&\delta(j,xy)=\delta(j,x)y+(-1)^{|x|}x\delta(j,y),\\
    &&\delta(k,xy)=\delta(k,x)y+(-1)^{|x|}x\delta(k,y).
\end{eqnarray*}
Then, $\delta(1,.),\delta(i,.)\in Der_0(\mathcal{H}^{a,b})$ and $\delta(j,.),\delta(k,.)\in Der_1(\mathcal{H}^{a,b})$. Therefore, by Proposition \ref{pr1} and Proposition \ref{pr2} the matrix of $\delta(1,.),\delta(i,.),\delta(j,.),\delta(k,.)$ has respectively the following form
\begin{eqnarray*}
   && M_1=\begin{pmatrix}
   0&0&0&0\\
   0&0&0&0\\
   0&0&0&- a\lambda^1\\
   0&0&\lambda^1&0\end{pmatrix}, \quad
   M_i=\begin{pmatrix}
   0&0&0&0\\
   0&0&0&0\\ 
   0&0&0&- a\lambda^i\\
   0&0&\lambda^i&0
   \end{pmatrix}\\
  && M_j=\begin{pmatrix}
  0&0&\mu^j&\nu^j\\
  0&0&0& 0 \\
  0&-b^{-1}\nu^j&0&0\\ 
  0&b^{-1}\mu^j&0&0
  \end{pmatrix}, \quad 
  M_k=\begin{pmatrix}
  0&0&\mu^k&\nu^k\\
  0&0&0& 0 \\ 
  0&-b^{-1}\nu^k&0&0\\ 
  0&b^{-1}\mu^k&0&0
  \end{pmatrix}
\end{eqnarray*}
From the equalities $\delta(1,j)=-\delta(j,1)$ we get 
\begin{equation}\label{equ1}
   a^1=0. 
\end{equation}
From the equalities $\delta(i,j)=-\delta(j,i)$ we get 
\begin{equation}\label{equ2}
   \lambda^i=-b^{-1}\mu^j, \quad \nu^j=0. 
\end{equation}
From the equalities $\delta(i,k)=-\delta(k,i)$ we deduce 
\begin{equation}\label{equ3}
   -a \lambda^i=b^{-1}\nu^k, \quad \mu^k=0. 
\end{equation}
From the equalities $\delta(j,k)=-(-1)^{|j||k|}\delta(k,j)$ we get 
\begin{equation}\label{equ4}
   \nu^j=\mu^k .
\end{equation}
By equations (\ref{equ1})-(\ref{equ4}) and setting, $a=a^i$ we have $\mu^j=-b \lambda$, $\nu^k=-ab \lambda$ and $\lambda^1=\nu^j=\mu^k=0$. Thus
\begin{eqnarray*}
   && M_1=\begin{pmatrix}0&0&0&0\\ 0&0&0&0\\ 0&0&0&0\\ 0&0&0&0\end{pmatrix}, \quad M_i=\begin{pmatrix}0&0&0&0\\ 0&0&0&0\\ 0&0&0&-a \lambda\\ 0&0&\lambda&0\end{pmatrix} \\
  && M_j=\begin{pmatrix}0&0&-b \lambda&0\\ 0&0&0& 0 \\ 0&0&0&0\\ 0&-\lambda&0&0\end{pmatrix}, \quad M_k=\begin{pmatrix}0&0&0&-ab \lambda\\ 0&0&0& 0 \\ 0&a \lambda&0&0\\ 0&0&0&0\end{pmatrix}
\end{eqnarray*}
Therefor, 
\begin{eqnarray*}
\delta(x,y)&=&x_2\delta(i,y)+x_3\delta(j,y)+x_4\delta(k,y)=x_2\widehat{M_i\overline{y}}+x_3\widehat{M_j\overline{y}}+x_4\widehat{M_k\overline{y}}\\
&=&x_2\widehat{(0,0,-a \lambda y_4,\lambda y_3)}+
x_3\widehat{(-b \lambda y_3,0,0,-\lambda y_2)}+
x_4\widehat{(-ab \lambda y_4,0,a \lambda y_2,0)}\\
&=&-ab\lambda(x_3y_3+x_4y_4)+a\lambda(x_4y_2-x_2y_4)j+\lambda(x_2y_3-x_3y_2)k.
\end{eqnarray*}
\end{proof}
We now prove that any symmetric  superderivations of degree $0$ for the generalized quaternion algebra is the zero map.
\begin{theorem}\label{thb2}
    Let $\delta$ be a symmetric super-biderivation of degree $0$ on $\mathcal{H}^{a,b}$. Then,    $$\delta(x,y)=0, \quad \forall x,y \in \mathcal{H}^{a,b}.$$
\end{theorem}
\begin{proof}
Let $\delta$ be a skew-super-biderivation of degree $0$ on  $(\mathcal{H}^{a,b})$, then for any $x,y \in \mathcal{H}^{a,b}$ we have 
\begin{eqnarray*}
    &&\delta(1,xy)=\delta(1,x)y+x\delta(1,y),\\
    &&\delta(i,xy)=\delta(i,x)y+x\delta(i,y),\\
    &&\delta(j,xy)=\delta(j,x)y+(-1)^{|x|}x\delta(j,y),\\
    &&\delta(k,xy)=\delta(k,x)y+(-1)^{|x|}x\delta(k,y).
\end{eqnarray*}
Then, $\delta(1,.),\delta(i,.)\in Der_1(\mathcal{H}^{a,b})$ and $\delta(j,.),\delta(k,.)\in Der_0(\mathcal{H}^{a,b})$. Therefore, by Proposition \ref{pr1} and Proposition \ref{pr2} the matrix of $\delta(1,.),\delta(i,.),\delta(j,.),\delta(k,.)$ has respectively the following form
\begin{eqnarray*}
&& M_1=\begin{pmatrix}0&0&0&0\\ 0&0&0&0\\ 0&0&0&-a \lambda^1\\ 0&0&\lambda^1&0\end{pmatrix}, \quad M_i=\begin{pmatrix}0&0&0&0\\ 0&0&0&0\\ 0&0&0&-a \lambda^i\\ 0&0&\lambda^i&0\end{pmatrix}\\
   && M_j=\begin{pmatrix}
  0&0&\mu^j&\nu^j\\
  0&0&0& 0 \\
  0&-b^{-1}\nu^j&0&0\\ 
  0&b^{-1}\nu^j&0&0
  \end{pmatrix}, \quad 
  M_k=\begin{pmatrix}
  0&0&\mu^k&\nu^k\\
  0&0&0& 0 \\ 
  0&-b^{-1}\nu^k&0&0\\ 
  0&b^{-1}\mu^k&0&0
  \end{pmatrix}  
\end{eqnarray*}
From the equalities $\delta(1,j)=\delta(j,1)$, we obtain
\begin{equation}\label{equ5}
    \lambda^1=0.
\end{equation}
From the equalities $\delta(i,j)=\delta(j,i)$, we have
\begin{equation}\label{equ6}
    \nu^j=0, \quad \lambda^i=b^{-1}\mu^j.
\end{equation}
From the equalities $\delta(i,k)=\delta(k,i)$, we have
\begin{equation}\label{equ7}
    \mu^k=0,\quad a \lambda^i=b^{-1}\nu^k.
\end{equation}
From the equalities $\delta(j,k)=-\delta(k,j)$, we have
\begin{equation}\label{equ8}
    \nu^j=-\mu^k.
\end{equation}
From the equalities $\delta(j,j)=-\delta(j,j)$, we get
\begin{equation}\label{equ9}
    \mu^j=-\mu^j.
\end{equation}
From the equalities $\delta(k,k)=-\delta(k,k)$, we have
\begin{equation}\label{equ8}
    \nu^k=-\nu^k.
\end{equation}
By comparing equations (\ref{equ5})-(\ref{equ10}) we deduce that
$$M_1=M_i=M_j=M_k=\begin{pmatrix}
  0&0&0&0\\
  0&0&0& 0 \\
  0&0&0&0\\ 
  0&0&0&0
  \end{pmatrix}.$$
This completes the proof.
\end{proof}
Combining Theorems \ref{thb1} and Theorem \ref{thb2} we have the following theorem
\begin{theorem}\label{thb3}
  Let $\delta$ be a super-biderivation of degree $0$ on $\mathcal{H}^{a,b}$. Then, there exists $\lambda \in \mathcal{R}$ such that 
    $$\delta(x,y)=-ab\lambda(x_3y_3+x_4y_4)+a\lambda(x_4y_2-x_2y_4)j+\lambda(x_2y_3-x_3y_2)k, \quad \forall x,y \in \mathcal{H}^{a,b}.$$  
\end{theorem}
\subsection{Characterization of $BDer_1(\mathcal{H}^{a,b})$}
We now present our second theorem on the structure of  superderivations of degree $1$ for the generalized quaternion algebra.\\
We start by proving that any skew super-biderivation of degree $1$ on the generalized quaternion algebra is the zero map.
\begin{theorem}\label{thb4}
    Let $\delta$ be a skew-super-biderivation of degree $1$ on $\mathcal{H}^{a,b}$. Then, 
    $$\delta(x,y)=0, \quad \forall x,y \in \mathcal{H}^{a,b}.$$
\end{theorem}
\begin{proof}
Let $\delta$ be a skew-super-biderivation of degree $1$ on  $(\mathcal{H}^{a,b})$, then for any $x,y \in \mathcal{H}^{a,b}$ we have 
\begin{eqnarray*}
    &&\delta(1,xy)=\delta(1,x)y+(-1)^{|x|}x\delta(1,y),\\
    &&\delta(i,xy)=\delta(i,x)y+(-1)^{|x|}x\delta(i,y),\\
    &&\delta(j,xy)=\delta(j,x)y+x\delta(j,y),\\
    &&\delta(k,xy)=\delta(k,x)y+x\delta(k,y).
\end{eqnarray*}
Then, $\delta(1,.),\delta(i,.)\in Der_1(\mathcal{H}^{a,b})$ and $\delta(j,.),\delta(k,.)\in Der_0(\mathcal{H}^{a,b})$. Therefore, by Proposition \ref{pr1} and Proposition \ref{pr2} the matrix of $\delta(1,.),\delta(i,.),\delta(j,.),\delta(k,.)$ has respectively the following form
\begin{eqnarray*}
   &&M_1=\begin{pmatrix}
  0&0&\mu^1&\nu^1\\
  0&0&0& 0 \\
  0&-b^{-1}\nu^1&0&0\\ 
  0&b^{-1}\mu^1&0&0
  \end{pmatrix}, \quad 
  M_i=\begin{pmatrix}
  0&0&\mu^i&\nu^i\\
  0&0&0& 0 \\ 
  0&-b^{-1}\nu^i&0&0\\ 
  0&b^{-1}\mu^i&0&0
  \end{pmatrix} \\
  && M_j=\begin{pmatrix}
   0&0&0&0\\
   0&0&0&0\\
   0&0&0&-a \lambda^j\\
   0&0&\lambda^j&0\end{pmatrix}, \quad
   M_k=\begin{pmatrix}
   0&0&0&0\\
   0&0&0&0\\ 
   0&0&0&-a \lambda^k\\
   0&0&\lambda^k&0
   \end{pmatrix}
\end{eqnarray*}
From the equalities $\delta(1,j)=-\delta(j,1)$ we get 
\begin{equation}\label{e1}
   \mu^1=\nu^1=0. 
\end{equation}
From the equalities $\delta(i,j)=-\delta(j,i)$ we get 
\begin{equation}\label{e2}
   \mu^i==0. 
\end{equation}
From the equalities $\delta(i,k)=-\delta(k,i)$ we get 
\begin{equation}\label{e3}
   \nu^i==0. 
\end{equation}
From the equalities $\delta(j,k)=\delta(k,j)$ we deduce 
\begin{equation}\label{e4}
   \lambda^j=\lambda^k=0. 
\end{equation}
By equations (\ref{e1})-(\ref{e4}) we deduce
\begin{eqnarray*}
    M_1=M_i=M_j=M_k=\begin{pmatrix}0&0&0&0\\ 0&0&0&0\\ 0&0&0&0\\ 0&0&0&0\end{pmatrix}.
\end{eqnarray*}
Therefor, 
$\delta(x,y)=x_1\widehat{M_1\overline{y}}+x_2\widehat{M_i\overline{y}}+x_3\widehat{M_j\overline{y}}+x_4\widehat{M_k\overline{y}}
=0$. This completes the proof.
\end{proof}
Now, we go to prove that any symmetric super-biderivation of degree $1$ on the generalized quaternion algebra is the zero map.
\begin{theorem}\label{thb5}
    Let $\delta$ be a symmetric super-biderivation of degree $0$ on $\mathcal{H}^{a,b}$. Then,    $$\delta(x,y)=0, \quad \forall x,y \in \mathcal{H}^{a,b}.$$
\end{theorem}
\begin{proof}
Let $\delta$ be a symmetric super-biderivation of degree $1$ on  $(\mathcal{H}^{a,b})$, then for any $x,y \in \mathcal{H}^{a,b}$ we have 
\begin{eqnarray*}
    &&\delta(1,xy)=\delta(1,x)y+(-1)^{|x|}x\delta(1,y),\\
    &&\delta(i,xy)=\delta(i,x)y+(-1)^{|x|}x\delta(i,y),\\
    &&\delta(j,xy)=\delta(j,x)y+x\delta(j,y),\\
    &&\delta(k,xy)=\delta(k,x)y+x\delta(k,y).
\end{eqnarray*}
Then, $\delta(1,.),\delta(i,.)\in Der_1(\mathcal{H}^{a,b})$ and $\delta(j,.),\delta(k,.)\in Der_0(\mathcal{H}^{a,b})$. Therefore, by Proposition \ref{pr1} and Proposition \ref{pr2} the matrix of $\delta(1,.),\delta(i,.),\delta(j,.),\delta(k,.)$ has respectively the following form
\begin{eqnarray*}
&& M_1=\begin{pmatrix}
  0&0&\mu^1&\nu^1\\
  0&0&0& 0 \\
  0&-b^{-1}\nu^1&0&0\\ 
  0&b^{-1}\mu^1&0&0
  \end{pmatrix}, \quad 
  M_i=\begin{pmatrix}
  0&0&\mu^i&\nu^i\\
  0&0&0& 0 \\ 
  0&-b^{-1}\nu^i&0&0\\ 
  0&b^{-1}\mu^i&0&0
  \end{pmatrix}\\
   && M_j=\begin{pmatrix}0&0&0&0\\ 0&0&0&0\\ 0&0&0&-a \lambda^j\\ 0&0&\lambda^j&0\end{pmatrix}, \quad M_k=\begin{pmatrix}0&0&0&0\\ 0&0&0&0\\ 0&0&0&-a \lambda^k\\ 0&0&\lambda^k&0\end{pmatrix}  
\end{eqnarray*}
From the equalities $\delta(1,i)=\delta(i,1)$, we obtain
\begin{equation}\label{e5}
    \mu^1=\nu^1=0.
\end{equation}
From the equalities $\delta(i,j)=\delta(j,i)$, we have
\begin{equation}\label{e6}
    \mu^i=0.
\end{equation}
From the equalities $\delta(i,k)=\delta(k,i)$, we have
\begin{equation}\label{e7}
    \nu^i=0.
\end{equation}
From the equalities $\delta(j,k)=-\delta(k,j)$, we have
\begin{equation}\label{e8}
    \lambda^j=\lambda^k=0.
\end{equation}
By equations (\ref{e5})-(\ref{e8}) we deduce that
$$M_1=M_i=M_j=M_k=\begin{pmatrix}0&0&0&0\\ 0&0&0&0\\ 0&0&0&0\\ 0&0&0&0\end{pmatrix}.$$
This completes the proof.
\end{proof}
Combining Theorems \ref{thb4} and \ref{thb5} we have the following theorem
\begin{theorem}
 Let $\mathcal{H}^{a,b}$ the quaternion algebra, then $$BDer_1(\mathcal{H}^{a,b})=\lbrace0\rbrace.$$   
\end{theorem}

where $D_{0}$ is a superderivation of degree 0 and $D_{1}$ is a superderivation of degree 1. Hence $[D]=[D_{0}]+[D_{1}]$. From Lemma 3.1, $D_{0}(1)=D_{0}(i)=0$ so we have $d_{11}=d_{12}=d_{13}=d_{14}=d_{21}=d_{22}=d_{23}=d_{24}=0.$ Using $D_{0}(j)=ak$ and $D_{0}(k)=-\alpha aj$ for some $a\in R$, we conclude that $d_{34}=a$ and $d_{43}=-\alpha a$. Therefore,
$$[D_{0}]=\begin{pmatrix}0&0&0&0\\ 0&0&0&0\\ 0&0&0&-\alpha a\\ 0&0&a&0\end{pmatrix}$$
Similarly, using Lemma 3.2, $D_{1}(1)=0$ and $D_{1}(j)=b$ for $b\in R$. Hence, in $[D_{1}]$ we see that $d_{11}=d_{12}=d_{13}=d_{14}=d_{32}=d_{33}=d_{34}=0$ and $d_{31}=b.$ Also, $D_{1}(k)=c$ implies $d_{41}=c$ and $d_{42}=d_{43}=d_{44}=0$. In addition, using $D_{1}(i)=-\beta^{-1}cj+\beta^{-1}bk$ we arrive at $d_{23}=-\beta^{-1}c$ and $d_{24}=\beta^{-1}b$. Thus,
$$[D_{1}]=\begin{pmatrix}0&0&b&c\\ 0&0&-\beta^{-1}c& \beta^{-1}b \\ 0&0&0&0\\ 0&0&0&0\end{pmatrix}$$
$$[D] = [D_0] + [D_1] = \begin{pmatrix}0&0&0&0\\ 0&0&0&0\\ 0&0&0&-\alpha a\\ 0&0&a&0\end{pmatrix} + \begin{pmatrix}0&0&b&c\\ 0&0&-\beta^{-1}c& \beta^{-1}b \\ 0&0&0&0\\ 0&0&0&0\end{pmatrix}$$
$$[D_{1}]=\begin{pmatrix}0&0&b&c\\ 0&0&0&0\\ 0&-\beta^{-1}c&0&0\\ 0&\beta^{-1}b&0&0\end{pmatrix} \quad \text{(as shown in the source 615)}$$
And the result holds:
$$[D] = [D_{0}] + [D_{1}] = \begin{pmatrix}0&0&b&c\\ 0&0&0&0\\ 0&-\beta^{-1}c&0&-\alpha a\\ 0&\beta^{-1}b&a&0\end{pmatrix}$$
\end{proof}

As a consequence of Theorem 3.4, we have the following result:

\textbf{Corollary 3.5.} The algebra $\mathrm{Der}_{s}(\mathcal{H}^{a,b})$ of superderivations of $\mathcal{H}^{a,b}$ is generated by the following matrices:
The algebra of superderivations $\mathrm{Der}_{s}(\mathcal{H}^{a,b})$ is generated by inner superderivations of the form given in Theorem 3.4; hence, the superderivation algebra coincides with the algebra of inner superderivations: $\mathrm{Der}_{s}(\mathcal{H}^{a,b})=\mathrm{Inn}_{s}(\mathcal{H}^{a,b})$. And the algebra of outer superderivations $\mathrm{Out}_{s}(\mathcal{H}^{a,b})$ vanishes for any generalized quaternion algebra $\mathcal{H}^{a,b}$.

\vspace{1cm}
\noindent \textbf{Acknowledgements:} The authors thank the referees for their valuable comments that contributed to a sensible improvement of the paper.

\section*{Declarations}
\noindent {\bf Conflict of interest statement.} Not applicable
\section*{Data availability}
The paper has no associated data.

\end{document}